\documentclass[11pt]{amsart}
\usepackage{amssymb,latexsym,comment,url,amsmath}

\usepackage[T1]{fontenc}
\usepackage{rotating,graphicx}

\usepackage{diagbox}

\usepackage{xcolor}

\usepackage{longtable}

\newcommand{\isom}{ \cong }

\newcommand{\PP}{{\mathbb P}}
\newcommand{\Q}{{\mathbb Q}}

\newcommand{\Orb}{\operatorname{Orb}}

\newcommand{\Z}{{\mathbb Z}}

\newcommand{\Magma}{{\sf MAGMA }}

\newcommand{\Mathematica}{{\sf Mathematica }}

\newenvironment{Proof}{\par\noindent{\sc Proof:}}%
{\hspace*{\fill}\nobreak$\Box$\par\medskip}
{\hspace*{\fill}\nobreak$\Box$\par\medskip}

\newtheorem{Proposition}{Proposition}[section]
\newtheorem{Theorem}[Proposition]{Theorem}

\newtheorem{Corollary}[Proposition]{Corollary}

\theoremstyle{definition}

\newtheorem{Remark}[Proposition]{Remark}
\newtheorem{Example}[Proposition]{Example}
\newtheorem{Conjecture}[Proposition]{Conjecture}

\renewcommand{\baselinestretch}{1.1}

\begin{document}

	\title[A Dynamical Analogue of a question of Fermat]{A Dynamical Analogue of a question of Fermat}

	\author[M. Sadek]%
	{Mohammad~Sadek}
	\address{Faculty of Engineering and Natural Sciences, Sabanc{\i} University, Tuzla, \.{I}stanbul, 34956 Turkey}
	\email{mmsadek@sabanciuniv.edu}

	\author[T. Yes\.{I}n]%
	{Tu\u{g}ba Yes\.{I}n}

	\email{tugbayesin@sabanciuniv.edu}

	\date{}

	\begin{abstract}
 Given a quadratic polynomial with rational coefficients, we investigate the existence of consecutive squares in the orbit of a rational point under the iteration of the polynomial. We display three different constructions of $1$-parameter quadratic polynomials with orbits containing three consecutive squares. 
		In addition, we show that there exists at least one polynomial of the form $x^2+c$ with a rational point whose orbit under this map contains four consecutive squares. This can be viewed as a dynamical analogue of a question of Fermat on rational squares in arithmetic progression.  
		Finally, assuming a standard conjecture on exact periods of periodic points of quadratic polynomials over the rational field, we give necessary and sufficient conditions under which the orbit of a periodic point contains only rational squares. 
	\end{abstract}
	
	\maketitle
\let\thefootnote\relax

	\section{Introduction}
	The existence of three consecutive squares in arithmetic progression is a phenomenon that can be seen in $\Q$. The rationals $1$, $5^2$, and $7^2$ provide such an example. In fact, one can parametrize all such rationals by observing that they satisfy the following equation
	$$x_2^2-x_1^2=x_3^2-x_2^2.$$ 
	This means that three rational numbers in arithmetic progression give rise to a rational point $(x_1:x_2:x_3)$ on the conic $C:x_1^2-2x_2^2+x_3^2=0$. Since the point $(1:1:1)\in C(\Q)$, it follows that $C(\Q)$ has infinitely many points. Moreover, one may parametrize these points as follows $(x_1:x_2:x_3)=(-p^2+2ps+s^2, p^2+s^2, p^2+2ps-s^2 )$ for some $ p,s\in \Q$.

	Fermat claimed that there does
	not exist an arithmetic progression of four squares over $\Q$. Euler, among others, proved this statement. One sees that the existence of such squares is equivalent to the existence of nontrivial rational points on the intersection of the following two quadric surfaces in $\PP^3_{\Q}$
	\begin{eqnarray*}
		x_1^2-2x_2^2+x_3^2&=&0\\
		x_2^2-2x_3^2+x_4^2&=&0.
	\end{eqnarray*}
	The latter intersection describes an elliptic curve $E$ for which $E(\Q)=\{(1 : \pm1 : \pm1 : \pm1)\}$. The points in $E(\Q)$ do not give rise to any non-constant rational squares in arithmetic progression. 
	
	In \cite{Xarles}, it was proved that a uniform upper bound exists on the number of
	squares in arithmetic progression over a given number field that depends
	only on the degree of the field. Moreover, the author proved that this bound is $5$ for
	quadratic fields. In \cite{Gonz}, the authors provide several criteria to identify the quadratic number fields over which there is a non-constant arithmetic progression of five squares.
	
	One may ask the aforementioned questions in a different setting, namely within the frame of  arithmetic dynamical systems. A dynamical system is a self-map $ f: S\longrightarrow S $ on a set $ S $ that allows iteration. The $ m$-${\textrm{th}} $ iterate of $ f $ is defined recursively by $ f^{0}(x)=x $ and $ f^{m}(x)=f(f^{m-1}(x)) $ when $ m\geq 1 $. The \textit{orbit} of a point $ P\in S $ under $ f $ is given by 
	\begin{equation*}
		\Orb_{f}(P)=\{f^{i}(P) : i=0,1,2\dots \}.
	\end{equation*}
	In case the map $ f $ is fixed, we write $\Orb(P)$. If $\Orb(P)$ is infinite, $P$ is called a wandering point; otherwise, $P$ is called a preperiodic point. A preperiodic point $P \in S$ is said to be {\em periodic} if there exists an integer $n > 0$ such that $f^n (P)= P$, where $n$ is called the period of $P$. If $n$ is the smallest such integer, we say that $P$ has exact period $n$. The orbit of a periodic point is called a {\em periodic orbit}.
	
	The question of the existence of $K$-rational squares in arithmetic progression of length $m$, $m\ge2$, over a number field $K$ can be reformulated using dynamical systems as follows. Can we find a linear polynomial $\ell(x)=x+c$, $c\in K^{\times}$, and $x_0\in K$ such that $\Orb_{\ell}(x_0)$ contains $m$ consecutive $K$-squares? In particular, is there an $x_0\in\Q$ such that $x_0,\ell(x_0),\ell^2(x_0),\ldots,\ell^{m-1}(x_0)$ are all in $K^2$?

	In this note, we are dealing with a higher degree dynamical analogue of Fermat's Squares Theorem. Namely, given a degree
	two polynomial $f(x) = x^{2} + Ax + B \in K[x]$ and a point $x_0\in K$,  how many consecutive squares can be there in the orbit $ \{x_0,f(x_0),f^2(x_0),  \ldots, f^n(x_0),\ldots \} $ of $x_0$?
	 It can be seen that for any irreducible quadratic map $f(x)\in K[x]$, the number of orbits under $f$ that contain at least three consecutive $K$-rational squares should be finite. This holds because each such square will give rise to a $K$-rational point on the hyperelliptic curve $C_m:y^2=f^m(x^2)$ for $m=0,1,2, \ldots$. When $m\ge 2$, the curve $C_m$ is of genus $>2$. By Faltings' theorem, see \cite{Faltings}, one then knows that the number of rational points on $C_m$, $m\ge 2$, must be finite.
	
	In this work, we give three different constructions of $1$-parameter polynomial maps of degree $2$ over $\Q$ and rational points that possess three different consecutive squares in their orbit under the iteration of these polynomials. In addition, unlike linear polynomial dynamical systems generated by polynomials of the form $x+c$, $c\in\Q^{\times}$, there exists at least one polynomial of the form $x^2+c$, $c\in \Q^{\times}$, and a point $x_0\in\Q$ such that $x_0,f(x_0), f^2(x_0)$ and $f^3(x_0)$ are all rational squares. 
	
	Finally, assuming a standard conjecture of Poonen on the exact period of periodic points of polynomial maps of degree $2$ over $\Q$, we introduce necessary and sufficient conditions under which polynomial maps of the form $x^2+ax+b\in\Q[x]$ possess periodic orbits containing only rational squares. 

 It is worth mentioning that a closely-related question was discussed in \cite{Cahn}. In fact, given an integer $r\ge 2$, the authors presented a complete classification of
rational functions defined over $K$ that possess an orbit containing infinitely many distinct $K$-rational $r$-th powers. 
	
	\subsection*{Acknowledgments} The authors would love to thank Mohamed Wafik for several suggestions. This work is supported by The Scientific and Technological Research Council of Turkey, T\"{U}B\.{I}TAK; research grant: ARDEB 1001/120F308. M. Sadek is partially funded by BAGEP Award of the Science Academy, Turkey.
	
	\section{Consecutive Three Squares}

  Let $K$ be a number field. Let $f\in K[x]$ and $x_0\in K$. 
  We say that $\Orb_{f}(x_0)$ contains {\em $m$-consecutive squares} if there is $y\in\Orb_f(x_0)$ such that
	\[y,f(y),\ldots,f^{m-1}(y)\] are all $K$-rational squares. We note that in the latter case  $\Orb_f(y)$ itself contains $m$-consecutive squares. Therefore, for the sake of simplicity, when we say that $\Orb_{f}(x_0)$ contains {\em $m$-consecutive squares} we mean
	\[x_0,f(x_0),\ldots,f^{m-1}(x_0)\] are all $K$-rational squares.

	We start with the following observation.
	\begin{Proposition}
		\label{prop1}
		Fix $a,b,c$ in a number field $K$. There are only finitely many $x_0\in K$ such that $\Orb_{f}(x_0)$, where $f(x)=ax^2+bx+c$, contains $3$ consecutive squares unless one of the following cases occurs.
 \begin{enumerate}
     \item $a=0$
     \item $b^2-4ac=0$
     \item $b=0$ and $c=\frac{-1}{a}$
     \item $b=4$ and $c=0$
     \item $1+2b\in \Q^{\times2}$ and $c=\frac{b^2-2b-2\pm2\sqrt{1+2b}}{4a}$
 \end{enumerate}
 Moreover, if $f(x)$ is an irreducible quadratic polynomial, then none of the cases above occurs, hence the finiteness of such $x_0$'s holds unconditionally.
 \end{Proposition}
	\begin{Proof}
		This follows immediately by observing that the existence of three consecutive squares can be expressed equivalently by
		\begin{eqnarray*}
			ax_0^4+bx_0^2+c=y^2,\qquad ay^4+by^2+c=z^2.
		\end{eqnarray*}
		This implies the existence of a rational point on the genus-$3$ curve
		\[C: z^2=a(ax_0^4+bx_0^2+c)^2+b(ax_0^4+bx_0^2+c)+c,\]
  
		By Faltings' Theorem, for fixed $K$-rational values $a,b,c$ such that the curve is smooth, the latter curve possesses only finitely many $K$-rational points. It remains to check the discriminant of the curve. Using $\Mathematica$, the discriminant is given by
  \[ \Delta=256 \, a^{15} c (1 + b + a c) (b^2 - 4 a c)^4 (-4 b^3 + b^4 + 16 a c + 16 a b c - 8 a b^2 c + 16 a^2 c^2)^2.\]
This gives the following cases for the curve not to be smooth:
\begin{enumerate}
    \item $a=0$, and for that case $f(x)$ is not a quadratic polynomial.
    \item $b^2-4ac=0$, in which case $f(x)=\frac{(b+2ax)^2}{4a}$. So, either $a$ is a square which gives that for any $x_0\in\Q^{\times2}$ the orbit will contain infinitely many consecutive squares; or $a$ is not a square, in which case for any $x_0\in\Q$, $f(x_0)$ is not a square.
    \item $c=0$ which gives rise to the curve
    \[C_1:Z_1^2:=\left(\frac{z}{x_0}\right)^2=(ax_0^2+b)(a^2x_0^4+abx_0^2+b).\]
    This is again a genus 2 curve and so Faltings' Theorem can still be applied, unless the discriminant of that curve given by $-64\, a^{15} (-4 + b)^2 b^8$ is zero. This gives that either $a=0$ covered in (i); $b=0$ implying that $b^2-4ac=0$ which is covered in (ii); or $b=4$ covered in (iv).
    \item $1+b+ac=0$ or $b=-1-ac$ gives rise to the curve
    \[C_2: Z_1^2:=\left(\frac{z}{x_0}\right)^2=(ax_0^2-ac-1)(a^2x_0^4-(a^2c+a)x_0^2+ac-1).\]
    The discriminant of the latter genus 2 curve is given by $64\, a^{15} (-1 + a c)^5 (1 + a c) (5 - 2 a c + a^2 c^2)^2$. If the curve is not smooth, then either $ac=1$ which means $b=-2$ and so $b^2-4ac=0$; $ac$ is a root of the irreducible polynomial $x^2-2x+5$, i.e, $ac\not\in\Q$; or $ac=-1$ which gives rise to (iii).
    \item Finally, the vanishing of the factor $-4 b^3 + b^4 + 16 a c + 16 a b c - 8 a b^2 c + 16 a^2 c^2$ in $\Delta$ yields that $ac$ is a root of the quadratic polynomial $16x^2+16(1+b-b^2)x+b^4-4b^3$ giving rise to the case (v).
\end{enumerate}
This concludes the proof. One can check easily that the aforementioned cases implies that $f(x)$ is reducible.
	\end{Proof}
	
	Two polynomials $f_1$ and $f_2$ are called $K$-linearly equivalent if there is a map $\ell(x)=ax+b\in K[x]$ such that $f_1=\ell\circ f_2\circ \ell^{-1}$. It is a simple exercise to see that any polynomial map of degree $2$ in $K[x]$ is $K$-linearly equivalent to map of the form $x^2+c$, $c\in K$. In what follows we focus on consecutive squares in orbits of points under maps of the form $f_c(x)=x^2+c$, $c\in\Q^{\times}$. 
	
	\begin{Theorem}
		For each $\beta\in\Q$, there are infinitely many rational numbers $\alpha,\gamma$, and $c$ such that $f_c(\alpha^2)=\beta^2$ and $f_c(\beta^2)=\gamma^2$.
		
		In particular,  one may choose
		\begin{eqnarray*}
			\alpha &=& \frac{\beta^2 (3 - 4 \beta^4)^2}{(1 + 8 \beta^2 + 4 \beta^4)^2},\\
			\gamma &=&\frac{\beta (-1 + 24 (\beta^2 + 3 \beta^4 + 4 \beta^6 + 2 \beta^8))}{(1 + 8 \beta^2 + 4 \beta^4)^2},\\
			c&=&\frac{\beta^2 - 49 \beta^4 + 400 \beta^6 + 2864 \beta^8 + 7264 \beta^{10} + 8864 \beta^{12} + 
				6400 \beta^{14} + 2816 \beta^{16} + 256 \beta^{18} - 256 \beta^{20}}{(1 + 8 \beta^2 + 4 \beta^4)^4}.
		\end{eqnarray*}
	\end{Theorem}
	\begin{Proof}
		Let $ \alpha \in \Q $ be such that $ f_c(\alpha)=\beta^{2} $ and  $ f_c(f_c(\alpha))=\gamma^{2} $ for some $ \gamma,\beta \in \Q, c\in\Q^{\times} $. This can be written as
		
		\begin{equation*}
			\alpha^{2}+c=\beta^{2}
		\end{equation*}
		\begin{equation*}
			\beta^{4}+c=\gamma^{2}.
		\end{equation*}
		Eliminating $c$, one has
		\begin{eqnarray}
			\alpha^{2}+\gamma^{2}=\beta^{4}+\beta^{2}.
		\end{eqnarray}
		For a fixed $ \beta $,  equation (1) defines a conic $C_{\beta}$ over $\Q(\beta)$ possessing a rational point $P_{\beta}$ defined by $  (\alpha, \gamma)=(\beta^{2}, \beta)\in C_{\beta}(\Q(\beta))$. Parameterizing the rational points $(\alpha,\gamma)\in C_{\beta}(\Q(\beta))$, using the point $P_{\beta}$, yields that 
		\begin{eqnarray*}
			\alpha=\frac{\beta(-2m-\beta+m^{2}\beta)}{1+m^{2}},\qquad \gamma=m(\alpha-\beta^2)+\beta=-\frac{\beta (-1 + 2 \beta m + m^2)}{1 + m^2}, \qquad\textrm{where }m\in \Q. 
		\end{eqnarray*} 
		Now, forcing $ \alpha $ to be a rational square, say $k^2$, we obtain the following quartic curve defined over $\Q(\beta)$ 
		\begin{equation*}
			H_{\beta}: k^2= \beta(\beta m^{4}-2 m^{3}z-2 mz^3-\beta z^4)
		\end{equation*}
		with the rational point $ (m:z:k)=(1 : 0 : \beta) $, hence $H_{\beta}$ is an elliptic curve over $\Q(\beta)$. The curve $H_{\beta}$ is $\Q$-birationally equivalent to the elliptic curve 
		\begin{equation*}
			E_{\beta}: y^2 = x^3 + (4\beta^4 + 4\beta^2)x.
		\end{equation*}
		We set $P_{\beta}=(1:0:\beta)$ and $ \phi:H_{\beta}\to E_{\beta} $ to be the birational isomorphism. One sees that $\phi(P_{\beta})=(1:2\beta^2+1:1)$ is of infinite order in $E_{\beta}(\Q(\beta))$ using \Magma, \cite{Magma}.  Proving the first part of the theorem.

		Now one has $ \phi^{-1}\left(2\phi(P_{\beta})\right)$ is given by
		\begin{equation*}
			\left( \left(-\frac{1}{2}\beta^4 - \frac{1}{8}\right)/\left(\beta^3 + \frac{1}{2}\beta\right) : \left(-\frac{1}{4}\beta^8 - \frac{1}{2}\beta^6 + \frac{1}{8}\beta^4 + \frac{3}{8}\beta^2 + \frac{3}{64}\right)/\left(\beta^5 + \beta^3 + \frac{1}{4}\beta\right) : 1 \right),
		\end{equation*}\\
		where the corresponding $ m$-coordinate on $H_{\beta}$ must be $ \left(-\frac{1}{2}\beta^4 - \frac{1}{8}\right)/\left(\beta^3 + \frac{1}{2}\beta\right) $.
		Consequently, one has the values given in the theorem.
	\end{Proof}
	
	\begin{Corollary}
		There are infinitely many $c\in \Q$ such that for some $x_0 \in \Q$, the orbit $ \Orb_{f_c}(x_0) $, where $ f_c(x)=x^2+c $, has three distinct consecutive squares.
	\end{Corollary}
	\begin{Example}
		Setting $ \beta =2 $, it can be seen that
		for the quadratic map $  f(x)=x^2+ 132583668/88529281 $ and $\alpha=122/97$, one has
		\begin{equation*}
			f_c(\alpha^2)=2^{2} \textrm{ and }f_c(4)=(39358/9409)^{2}.
		\end{equation*} 
	\end{Example}

	\begin{Theorem}
	Let $a\in \Q$. There exist infinitely many $\delta,\gamma,b\in\Q$ such that for the map $f(x)=x^2+ax+b$, one has $f(\delta^2)=a^2$ and $f(a^2)=\gamma^2$. In particular, one can choose
	\begin{equation*}
			b=\frac{a^2 - a^3 - 69 a^4 - 196 a^5 + 314 a^6 + 2226 a^7 + 7622 a^8 + 
				15308 a^9 + 25285 a^{10} + 30279 a^{11} + 31599 a^{12}}{(1 + a (2 + a (9 + 4 a (1 + a))))^4} 
		\end{equation*}
		\begin{equation*}
			+ \frac{24864 a^{13} + 16624 a^{14} + 6496 a^{15} + 160 a^{16} - 3072 a^{17} - 
				2560 a^{18} - 1280 a^{19} - 256 a^{20}}{(1 + a (2 + a (9 + 4 a (1 + a))))^4}
		\end{equation*}
	and
	\begin{equation*}
			\delta=\frac{a(1+a)(-3+a+4a^3)}{(1+a(2+a(9+4a(1+a))))}, \ \ \ \ \gamma= \frac{a (-1 + a^2 (5 + 4 a (1 + a)) (6 + a (8 + 3 a (5 + 4 a (1 + a)))))}{(1 + a (2 + a (9 + 4 a (1 + a))))^2}.
		\end{equation*}
	
		It follows that there exist infinitely many polynomials $ f(x)=x^2+ax+b \in \Q[x] $ such that $ \Orb_{f}(x) $ contains three distinct consecutive squares for some $ x\in \Q $. 
	\end{Theorem}
	
	\begin{Proof}
		Let $ \alpha \in \Q $ and assume $ f(\alpha)= \alpha^{2}+a \alpha+b=\beta^{2} $ and $ f(f(\alpha))=\beta^{4}+a\beta^{2}+b=\gamma^{2} $ for some $ \beta,\gamma \in \Q $. By eliminating $b$, we have 
		
		\begin{equation*}
			\alpha^2+a\alpha-\beta^2=\beta^4+a\beta^2-\gamma^2.
		\end{equation*}
		
		One observes that setting $ \beta=a $, the equation above describes a conic $ C_{a} :\alpha^2+\gamma^2+a\alpha=a^4+a^3+a^2 $ over $ \Q(a) $ possessing a rational point $ P_{a} $ defined by $ (\alpha,\gamma)=(a^2,a)\in C_{a}(\Q(a)) $. We parameterize the rational points $ (\alpha,\gamma)\in C_{a}(\Q(a)) $ using the point $ P_{a} $ as follows
		
		\begin{equation*}
			\alpha=\frac{a(-1-a-2m+am^2)}{1+m^2},\ \ \ \ \gamma=m(\alpha -a^2)+a=-\frac{a (-1 + m + 2 a m + m^2)}{1 + m^2},\quad m\in\Q.
		\end{equation*}
		Now, forcing $ \alpha $ to be a rational square, say $ k^2 $, we obtain the following quartic curve defined over $ \Q(a) $
		
		\begin{equation*}
			H_{a}: k^2= a^2m^4-2am^3-am^2-2am-a^2-a
		\end{equation*}
		with a rational point $ (m:k:z)=(1,a,0)$. Therefore, $ H_{a} $ is an elliptic curve over $ \Q(a) $ and it is $ \Q $-birationally equivalent to the elliptic curve $E_a$ defined by the Weierstrass equation
		
		\begin{equation*}
			E_{a}: y^2-\frac{2}{a}xy-\frac{4a^2+2a+2}{a^3}y=x^3+\frac{2a+2}{a^2}x^2+\frac{4a^4+4a^3+a^2+2a+1}{a^4}x
		\end{equation*}
		One sees that the image of the point $Q_a=(1:a:0)$ in $E_a$ under the birational isomorphism $\psi:H_a\to E_a$  is $((0 : \frac{4a^2 + 2a + 2}{a^3} : 1))$ which is of infinite order, \Magma \cite{Magma}.
		Now the $ m$-coordinate of the rational point $\psi^{-1}(2\psi(Q_{a}))$ in $H_{a}$ is given by 
		
		\begin{equation*}
			 \left( -\frac{1}{2}a^4-\frac{1}{2}a^3-\frac{1}{8}a^2-\frac{1}{4}a-\frac{1}{8}     \right)  / \left(   a^3+\frac{1}{2}a^2+\frac{1}{2}a                               \right).
		\end{equation*}
		With the latter $ m $-coordinate, we get the values for $b$, $\delta$ and $\gamma$ as in the statement of the theorem.
	\end{Proof}
	
The following theorem also describes an explicit construction of three consecutive squares in the orbit of polynomials of the form $x^2+ax-a$.
	
	\begin{Theorem}
	Let $\alpha\in \Q$ and $a=-\frac{(-1 + \alpha) \alpha^2 (-9 + \alpha^2)}{(4 + \alpha - \alpha^2)^2}$. For the polynomial $f(x)=x^2+ax-a$, one has 
	$$\Orb_f(\alpha)=\left\{\alpha, \left(\frac{\alpha^2 -5\alpha}{\alpha^2 - \alpha - 4}\right)^2,  \left(\frac{3\alpha^5 - 13\alpha^4 + 13\alpha^3 - 15\alpha^2 + 12\alpha}{(\alpha^4 - 2\alpha^3 - 7\alpha^2 + 8\alpha + 16) (\alpha-1)}\right)^2\right\}.$$
		In particular, for any rational number $x_0\in\Q$, there exists an $a\in\Q$ such that the polynomial $ f(x)=x^2+ax-a $ satisfies $f(x_0^2)$ and $f^2(x_0^2)$ are rational squares.
	\end{Theorem}
	
	\begin{Proof}
		Let $ \alpha \in \Q $ be such that $ f(\alpha)= \alpha^{2}+a \alpha-a=\beta^{2} $ and $ f(f(\alpha))=\beta^{4}+a\beta^{2}-a=\gamma^{2} $ for some $ \beta,\gamma \in \Q $. One obtains
		\begin{equation*}
			a=\frac{\beta^{2}-\alpha^{2}}{\alpha-1}=\frac{\gamma^{2}-\beta^{4}}{\beta^{2}-1}.
		\end{equation*}
		This gives a certain level of confidence. 
		\begin{equation*}
			\alpha\beta^{4}+(-1-\alpha^{2})\beta^{2}+\alpha^{2}=(\alpha-1)\gamma^{2}
		\end{equation*}\\
		which defines the following quartic curve over $\Q(\alpha)$
		\begin{equation*}
			C_{\alpha}: 	\alpha(\alpha-1)\beta^{4}+(\alpha-1)(-1-\alpha^{2})\beta^{2}+(\alpha-1)\alpha^{2}=\theta^{2}
		\end{equation*}
		where $ \theta=(\alpha-1)\gamma $, with a rational point $ T_{\alpha}=(1:\alpha-1:1) $.
		
	There is a birational isomorphism $\psi:C_{\alpha}\to E_{\alpha}$ where $ E_{\alpha} $ is an elliptic curve described by the following Weierstrass equation over $\Q(\alpha)$
		
		\begin{equation*}
			E_{\alpha}: y^2 + (2\alpha + 2)xy + \frac{2\alpha^4 - 6\alpha^3 - 2\alpha^2 - 2\alpha}{\alpha - 1}y = x^3 + \frac{2\alpha^3 - 8\alpha^2 - 2\alpha}{\alpha - 1}x^2 + \frac{\alpha^6 - 8\alpha^5 + 14\alpha^4 + 4\alpha^3 + 5\alpha^2}{\alpha^2 - 2\alpha + 1}x
		\end{equation*}\\
		
		Then $ R_{\alpha}:=\psi(T_{\alpha})=(0 : (-2\alpha^4 + 6\alpha^3 + 2\alpha^2 + 2\alpha)/(\alpha - 1) : 1) $ is a point of infinite order in $E_{\alpha}(\Q(\alpha))$.
Moreover,		
		\begin{equation*}
			 \psi^{-1}(2R_{\alpha})=\left( \frac{\alpha^2 -5\alpha}{\alpha^2 - \alpha - 4}: \frac{3\alpha^5 - 13\alpha^4 + 13\alpha^3 - 15\alpha^2 + 12\alpha}{\alpha^4 - 2\alpha^3 - 7\alpha^2 + 8\alpha + 16} : 1 \right).
		\end{equation*}
		Now the $\beta$-coordinate of the latter rational point gives rise to the $a$-value and the corresponding orbit in the statement of the theorem. 
	\end{Proof}
	
	\section{Consecutive four squares}

	Let $K$ be a number field. Let $f(x)=x^2+c\in K[x]$ and $x_0\in K$. If one wants to force $x_0^2,f_c(x_0^2), f_c^2(x_0^2)$ and $f_c^3(x_0^2)$ to be all $K$-rationals, then this can be written as
	\begin{eqnarray}
		\label{eq}
		x_0^4+c=y^2,\qquad y^4+c=z^2, \qquad z^4+c=w^2. 
	\end{eqnarray}
	Equivalently, the existence of four consecutive squares in the orbit of a rational point under $f_c$ is equivalent to the existence of a rational point $(x_0,y,z,w)$ on the surface $\mathcal S$ defined by \begin{eqnarray}
	\label{eq1}
		z^2+x^4=y^2+y^4, \qquad w^2+y^4=z^2+z^4.
	\end{eqnarray}
	
	\begin{Proposition}
		Let $c\in\Q$ be such that $\Orb_{f_c}(x_0)$ contains four consecutive squares, i.e, $f^i(x_0)$, $i=0,1,2,3$, are all rational squares. Then $x_0\ne 0$.
	\end{Proposition}
	\begin{Proof}
		It can be seen that by eliminating $y$ and $z$ in (\ref{eq}), one obtains that 
		\begin{eqnarray*}w^2&=&f_c^2(x_0^2)=((x_0^4+c)^2+c)^2+c\\
			&=&x_0^{16}+4cx_0^{12}+2c(1+3c)x_0^8+4c^2(1+c)x_0^4+c(1+c+2c^2+c^3).  
		\end{eqnarray*}
		If $x_0=0$, then this means that $w^2=c(1+c+2c^2+c^3)$ which describes an elliptic curve $E$ over $\Q$, whose Mordell-Weil group $E(\Q)\isom\Z/3\Z$ corresponding to the point $(0,0)$ and the two points at infinity. None of these points gives rise to non-trivial four consecutive squares. 
	\end{Proof}
	
	\begin{Theorem}
	
		 There exists a polynomial $ f(x)=x^2+c \in \Q[x] $ such that there are four distinct consecutive squares in $ \Orb_{f}(x_0^2) $ for some $x_0\in \Q$ if and only if there exist rational solutions $ p,q,r \in \Q $ to the polynomial equation $ M(p,q,r)=x_0^4 $ where 
   
		{\footnotesize $ M(p,q,r)=
			-2048 p^7 q - 1536 p^6 q^2 - 768 p^5 q^3 + 128 p^4 q^4 + 192 p^3 q^5 + 96 p^2 q^6 + 16 p q^7 + q^8 + 4096 p^6 q r + 
			1280 p^5 q^2 r - 960 p^3 q^4 r - 448 p^2 q^5 r - 160 p q^6 r - 
			24 q^7 r - 512 p^6 r^2 - 1280 p^5 q r^2 + 1792 p^4 q^2 r^2 + 
			1664 p^3 q^3 r^2 + 1152 p^2 q^4 r^2 + 368 p q^5 r^2 + 76 q^6 r^2 + 
			768 p^5 r^3 - 2048 p^4 q r^3 - 1920 p^3 q^2 r^3 - 1280 p^2 q^3 r^3 - 
			384 p q^4 r^3 - 72 q^5 r^3 + 384 p^4 r^4 + 1728 p^3 q r^4 + 
			480 p^2 q^2 r^4 + 240 p q^3 r^4 - 10 q^4 r^4 - 704 p^3 r^5 - 
			64 p^2 q r^5 - 32 p q^2 r^5 + 24 q^3 r^5 + 64 p^2 r^6 - 
			112 p q r^6 + 44 q^2 r^6 + 64 p r^7 - 56 q r^7 + 17 r^8$.}
	\end{Theorem}
	\begin{Proof}
		In (\ref{eq1}), we set $ x_0^2=X $, \ $y^2=Y $ and $ z^2=Z$.
		Then we have the following equations
		\begin{eqnarray} \label{eqq}
			Y+Y^2=Z+X^2,
		\end{eqnarray}
		\begin{eqnarray} \label{eqqq}
			Z+Z^2=w^2+Y^2.
		\end{eqnarray}
		One may homogenize equation (\ref{eqqq}) and complete the square so that the equation may be written as
		\begin{equation*}
			\gamma^2=w^2+Y^2+\mu^2,\qquad\textrm{where }
		  \mu=\frac{T}{2}, \, \gamma=Z+\mu.
		  \end{equation*}
		Therefore, one may obtain the following parameterization.
		\begin{equation*}
			\gamma=s^2+t^2+u^2, \ \ w=2su,\ \ Y=s^2+t^2-u^2, \ \ \mu=2tu.
		\end{equation*}
		Since $ Z=\gamma-\mu $ and $ T=2\mu $, we have $ z^2=-2tu+s^2+t^2+u^2$ and $ T=4tu $. Also $ Y=y^2=s^2+t^2-u^2 $ yielding the following parametrization for $ t,s,u,y $:
		\begin{eqnarray*}
			s=4p^2+2qp-2pr-q^2+r^2&,&   t=4qp+q^2-2qr+r^2,\\
			u=4p^2+2qp-2pr+q^2-r^2&,&  y=4pr+2qr-q^2-r^2.
		\end{eqnarray*}

		It follows that $(\ref{eqq})$ in homogeneous form, $YT+Y^2=ZT+X^2$, may be written as
		\begin{equation*}
			x_0^4=y^2T+y^4-Tz^2=M(p,q,r)
		\end{equation*}\\
		where $ M(p,q,r)$ is given as in the statement of the theorem.
	\end{Proof}

	\begin{Theorem}
	\label{thm:1}
		There exists at least one polynomial of the form $ f(x)=x^2+c \in \Q[x]$ and $x_0\in \Q$ such that 
		$ \Orb_{f}(x_0) $ has four distinct consecutive squares. Namely, $ c=5103/4096 $ and
		$$\Orb_f((3/8)^2)= \{(3/8)^2, (9/8)^2, (27/16)^2, (783/256)^2,\ldots\}.$$
	\end{Theorem}
	
	\begin{Proof}
		Fixing $ y $ in equations (\ref{eq1}) and setting $ X=x_0^2 $, one obtains	
		\begin{equation*}
			(y^2+y^4)T^2=z^2+X^2.
		\end{equation*}
		The latter equation gives rise to the following parametrization
		\begin{equation*}
			z=-yp^2+2y^2ps+ys^2, \ \ \ \ X=y^2+p^2+2yps-y^2s^2, \ \ \ \ T=p^2+s^2.
		\end{equation*}\\
		
		Let $ X_{1}=\frac{X}{T} $ and $ z_{1}=\frac{z}{T} $. Then we have the following. 
		\begin{equation*}
			X_{1}=\frac{z(2ps+p^2z-s^2z)}{p^2+s^2}=\square,
		\end{equation*}
		
		\begin{equation*}
			z_{1}^{2}+z_{1}^{4}-y^4=\frac{y^4(-p^2+s^2+2psy)^4+y^2(-p^2+s^2+2psy)^2(p^2+s^2)^2-y^4(p^2+s^2)^4}{(p^2+s^2)^4}=\square.
		\end{equation*}\\
		
		Searching for rational solutions to the system above using \Magma, \cite{Magma}, yields the polynomial $f(x)$ together with the mentioned orbit.  
	\end{Proof}
	
\begin{Remark}
In Theorem \ref{thm:1}, we were able to find a rational point on the variety $\mathcal S$ defined in (\ref{eq1}). This variety contains the subvariety (up to sign) $x=y=z=w$ with infinitely many rational points that give rise to no nontrivial four distinct consecutive squares. We suspect that there are likely only finitely many other nontrivial rational points, and perhaps the rational point we found might be the only one.     
\end{Remark}
As for polynomials $f(x)$ with $d=\deg f>2$, the existence of a rational square $\alpha^2$ such that $f\left(\alpha^2\right)$ is rational itself, implies the existence of a rational point on a curve of genus $\lfloor 2d-1\rfloor/2>1$, on which there are only finitely many rational points.  Therefore, finding a rational point whose orbit under $f$ contains three consecutive squares is quite improbable.  

\section{Finite orbits consisting of squares}

As mentioned before, any quadratic polynomial map $f(x) = Ax^2 + Bx + C \in K[x]$ is linearly conjugate over $K$ to a map of the form $x^2+ c$ for some $c \in K$. If $K$ is chosen to be the rational field $\Q$, a complete classification of quadratic polynomial maps with periodic points of periods $1, 2$, or $3$ was given in \cite{WaldeRusso}. We recall that the orbit of a periodic point is called a {\em periodic orbit}. The following can be found for example as \cite[Theorem 1]{Poonen}.
\begin{Proposition}
\label{prop:poonen}
Let $f(x)=x^2 +c$ with $c\in \Q$. Then
\begin{itemize}
\item[1)] $f(x)$ has a rational point of period $1$, i.e., a rational fixed point, if and only if $c = 1/4 -\rho^2$ for some $\rho\in  \Q$. In this case, there are exactly two, $1/2 +\rho$ and $1/2 -\rho$, unless $\rho = 0$, in which case they coincide.
\item[2)] $f(x)$ has a rational point of period $2$ if and only if $c = -3/4 -\sigma^2$  for some $\sigma\in\Q$, $\sigma\ne 0$. In this case, there are exactly two, $-1/2 +\sigma$ and $-1/2 -\sigma$  (and these form a $2$-cycle).
\item[3)] $f(x)$ has a rational point of period $3$ if and only if
$$c = -\frac{\tau^6  + 2\tau^ 5 + 4\tau^4 + 8\tau^3 + 9\tau^2 + 4\tau + 1}{
4\tau^2(\tau +1)^2}$$
for some $\tau  \in \Q$, $\tau\ne  -1, 0$. In this case, there are exactly three,
\begin{eqnarray*}
x_1 = \frac{\tau^3 + 2\tau^2 +\tau +1}{ 2\tau(\tau +1)},\quad
x_2 = \frac{\tau^3 -\tau -1}{ 2\tau(\tau +1)},\quad
x_3 = -\frac{ \tau^ 3 + 2\tau^2 + 3\tau + 1}{ 2\tau(\tau +1)}
\end{eqnarray*}
and these are cyclically permuted by $f(x)$.
\end{itemize}
\end{Proposition}   
The following conjecture can be found in \cite{Poonen}.
\begin{Conjecture}\label{Con:Poonen}
If $N \ge 4$, then there is no quadratic polynomial $f (x) \in\Q[x]$ with a rational point of exact period $N$ .
\end{Conjecture}
The conjecture has been proved for $N=4$, \cite{Morton}, for $N=5$, \cite{fps}, and conditionally on Birch-Swinnerton-Dyer Conjecture for $N=6$, \cite{Stoll}. Many results have been obtained on the size of the intersection of orbits of two degree-$2$
rational maps assuming that the latter conjecture holds true, \cite{Burcu,Hindes}. Although proving the uniform boundedness of the number of preperiodic points of rational maps of a fixed degree is currently far from our reach, some uniform bounds were given for certain polynomial and rational maps in \cite{Ingram,Sadek,Sadek1}.  

		Assuming Conjecture \ref{Con:Poonen} holds, one notices that if $ f(x)=x^2+c \in \Q[x]$ is such that $x_0\in \Q$ is a periodic point of $f(x)$, then for $x_0$ to be a rational square of period $1$ one has either $1/2+\rho$ or $1/2-\rho$ is a rational square with $c=1/4-\rho^2$. Similarly, one sees easily that $x_0$ cannot be a point of period $2$ whose orbit contains only rational squares since otherwise both $-1/2+\sigma$ and $-1/2-\sigma$ are rational squares for some $\sigma\in\Q$. Finally, for $x_0$ to be a periodic point of period $3$ for which $\Orb_f(x_0)$ contains only rational squares, one must have in Proposition \ref{prop:poonen} that $x_1=r_1^2, \ x_2=r_2^2,\ x_3=r_3^2$ where $r_i\in\Q$, $i=1,2,3$. The latter is a singular curve of genus $17$ with only two singularities $(\tau,r_1,r_2,r_3)=(-1, 0, 0, 0), (0, 0, 0, 0)$. Again, by Faltings' Theorem, \cite{Faltings}, there are only finitely many rational points on the latter curve. Therefore, one investigates the possibility of having infinitely many polynomials of the form $x^2+ax+b,\ a\ne 0,$ with rational periodic points whose orbits are of length at least $2$ and contain only rational squares.

  One notices that if $x_0$ is a rational periodic point of the map $x^2+ax+b$, then $x_0+a/2$ is a periodic point of the map $x^2+c$ where $c=b-a^2/4+a/2$.

\begin{Theorem} 
The polynomial map $f(x)=x^2+ax+b\in\Q[x]$ has a periodic orbit of length $2$ whose elements are rational squares if and only if $a=-1-m^2-k^2$ and $b=m^2+k^2+m^2k^2$ for some $m,k\in\Q$. In this case, one has $f(m^2)=k^2$ and $f(k^2)=m^2$.
\end{Theorem}
\begin{Proof}
    That the polynomial $f(x)=x^2+ax+b$ with $a=-1-m^2-k^2$ and $b=m^2+k^2+m^2k^2$, $m,k\in\Q$, has such a periodic orbit is a direct calculation. 

    Now, if $f(k^2)=m^2$ and $f(m^2)=k^2$ for some $m,k\in\Q$, then one knows that $g(m^2+a/2)=k^2+a/2$ and $g(k^2+a/2)=m^2+a/2$ where $g(x)=x^2+b-a^2/4+a/2$. 
 This yields that  
 \begin{equation*}
    k^2+\frac{a}{2}=-\frac{1}{2}-\sigma \ \ \text{and} \ \  m^2+\frac{a}{2}=-\frac{1}{2}+\sigma, \qquad \textrm{for some }\sigma\in\Q,
 \end{equation*}
see Proposition \ref{Con:Poonen}.
It follows that $a=-1-m^2-k^2$ and $b=m^2+k^2+m^2k^2$. 
\end{Proof}

One sees that $\Orb_{f}(4)=\{4,\frac{1}{4}\}$ where $f(x)=x^2-\frac{21}{4}x+\frac{21}{4}$.

\begin{Theorem}\label{thm} Let $m,n,r\in\Q$ be distinct. There exists a polynomial map $f(x)=x^2+ax+b\in\Q[x]$ such that $f(m^2)=n^2$, $f(n^2)=r^2$, and $f(r^2)=m^2$ if and only if $$m^4 (1 - n^2 + r^2) + 
 m^2 (-n^2 + n^4 - r^2 (1 + r^2))+r^4 - n^4 (-1 + r^2) + n^2 r^2 (-1 + r^2)=0.$$ 
In this case, the polynomial $f(x)$ is determined by 
$$a=\frac{-m^6 + m^2 n^4 - n^6 + m^4 r^2 + 
 n^2 r^4 - r^6}{(m^2 - n^2) (m^2 - r^2) (n^2 - r^2)},\qquad b= \frac{
m^6 n^2 - m^4 n^4 + n^6 r^2 - m^4 r^4 - n^4 r^4 + 
 m^2 r^6}{(-m^2 + n^2) (n^2 - r^2) (-m^2 + r^2)}.$$
\end{Theorem}
	\begin{Proof}
 One needs to solve the following system of linear equations in $d,a,b$
 \begin{eqnarray*}
 dm^4+am^2+b=n^2,\quad dn^4+an^2+b=r^2,\quad dr^4+ar^2+b=m^2
 \end{eqnarray*}
 to get the expressions for $a$ and $b$ as in the statement, whereas $d=(m^4 - m^2 n^2 + n^4 - m^2 r^2 - 
 n^2 r^2 + r^4)/((m^2 - n^2) (m^2 - r^2) (n^2 - r^2)).$ The statement now holds once we force the polynomial $f(x)$ to be monic by setting $d=1$.
	\end{Proof}
One remarks that each of the triples $m,n,k$ satisfying the identity in Theorem \ref{thm} gives rise to a rational solution to the system of equations   	
	\begin{eqnarray*}
m^2+a/2 = \frac{\tau^3 + 2\tau^2 +\tau +1}{ 2\tau(\tau +1)},\quad
n^2+a/2 = \frac{\tau^3 -\tau -1}{ 2\tau(\tau +1)},\quad
r^2+a/2 = -\frac{ \tau^ 3 + 2\tau^2 + 3\tau + 1}{ 2\tau(\tau +1)}
\end{eqnarray*}
for some rational value of $\tau\in\Q\setminus\{-1,0\}$. \\

As examples, one sees that the following polynomial maps have periodic orbits of length $3$ that contain only rational squares.
\begin{eqnarray*}
    f_1(x)=&x^2-\frac{29}{8} \ x+\frac{841}{256},\quad \Orb_{f_1}((7/4)^2)&=\{(7/4)^2,(5/4)^2,(1/4)^2\}, \quad \tau=-1/2,\\ 
    f_2(x)=&x^2-\frac{301}{72}\ x+\frac{90601}{20736},\quad\Orb_{f_2}((23/12)^2)&=\{(23/12)^2,(19/12)^2,(5/12)^2\},\quad\tau=2,\\
    f_3(x)=&x^2-\frac{421}{72} \ x+\frac{177241}{20736},\quad \Orb_{f_3}((25/12)^2)&=\{(25/12)^2,(17/12)^2,(11/12)^2\},\quad \tau=1/2,\\
    f_4(x)=&x^2-\frac{1849}{288} \ x+\frac{3418801}{331776},\quad \Orb_{f_4}((55/24)^2)&=\{(55/24)^2,(49/24)^2,(23/24)^2\},\quad \tau=3,\\
   f_5(x)=&x^2-\frac{74333}{4356} \ x+\frac{211660729}{4743684},\quad \Orb_{f_5}((115/66)^2)&=\{(115/66)^2,(47/33)^2,(124/33)^2\},\quad \tau=-12. \\
		\end{eqnarray*}

\end{document}